\documentclass[a4paper,oneside,12 pt]{article}
 
\usepackage[T1]{fontenc}
\usepackage[utf8]{inputenc}
\usepackage{amsthm}
\usepackage{amsmath}
\usepackage{amssymb}
\usepackage{mathrsfs}
\usepackage{graphicx}
\usepackage{enumitem} 
\usepackage{filecontents}
\usepackage{natbib}
\allowdisplaybreaks

\title{Weak boundedness of $\CZ$ operators on noncommutative $L_1$-spaces}
\author{Léonard Cadilhac}

\newtheorem{definition}{Definition}[section]
\newtheorem{theorem}[definition]{Theorem}
\newtheorem{property}[definition]{Proposition}
\newtheorem{lemma}[definition]{Lemma}

\theoremstyle{definition}
\newtheorem{remark}[definition]{Remark}

\newtheorem{corollary}[definition]{Corollary}

\newcommand{\R}{\mathbb{R}}
\newcommand{\Z}{\mathbb{Z}}
\newcommand{\Rn}{\mathbb{R}^n}

\newcommand{\Fn}{\mathbb{F}}
\newcommand{\Nb}{\mathbb{N}}
\newcommand{\Sb}{\mathbb{S}}
\newcommand{\Cb}{\mathbb{C}}

\newcommand{\A}{\mathcal{A}}
\newcommand{\B}{\mathcal{B}}
\newcommand{\Q}{\mathcal{Q}}

\newcommand{\M}{\mathcal{M}}
\newcommand{\N}{\mathcal{N}}
\newcommand{\E}{\mathcal{E}}

\newcommand\normi[1]{\big\| #1 \big\|_{\mathcal{B}(L_2(\N))}}
\newcommand\Norm[1]{\big\| #1 \big\|}
\newcommand\md[1]{\left| #1 \right|}

\newcommand\NormMt[1]{\big\| #1 \big\|_{\Mt}}

\newcommand{\les}{\lesssim}

\newcommand{\e}{\varepsilon}

\newcommand{\I}{\mathbf{1}}

\newcommand{\otimest}{\overline{\otimes}}

\newcommand{\Mt}{\widetilde{\mathcal{M}}}
\newcommand{\Tt}{\widetilde{T}}
\newcommand{\kt}{\widetilde{k}}
\newcommand{\Nt}{\widetilde{\mathcal{N}}}

\newcommand{\Zb}{\overline{\mathbb{Z}}}

\newcommand{\phih}{\widehat{\phi}}

\newcommand\supb[1]{\underset{#1}{\sup}}
\newcommand\Sum[2]{\sum\limits_{#1}^{#2}}
\newcommand\Int[2]{\int\limits_{#1}^{#2}}
\newcommand\Lim[1]{\underset{#1}{\lim}}

\newcommand\RC[2]{\text{RC}_{#1}(#2)}

\newcommand{\CZ}{\text{Calder\'on-Zygmund}}

\newcommand\addtag{\refstepcounter{equation}\tag{\theequation}}

\begin{document}

\maketitle

\begin{abstract}
In 2008, J. Parcet showed the $(1,1)$ weak-boundedness of $\CZ$ operators acting on functions taking values in a von Neumann algebra. We propose a simplified version of his proof using the same tools: Cuculescu's projections and a pseudo-localisation theorem. This will unable us to recover the $L_p$-boundedness of $\CZ$ operators with Hilbert valued kernels acting on operator valued functions for $1 < p < \infty$ and an $L_p$-pseudo-localisation result of P. Hytönen. 
\end{abstract}

\section{Introduction}

Before properly introducing the topic of this paper, it should be said that our main purpose is to offer a simplified version of the argument presented by J. Parcet in \cite{czoperator}. Consequently, its interest lies in the shorter and clearer proof it provides. We believe it can be used to expand the reach of the theorem as shown by the slight improvement we make to its hypotheses and an application we present.

Our main result belongs to the now well developped theory of singular integrals. The latter was initiated in the 1950's by \text{Calder\'on} and Zygmund who found a very useful sufficient condition for a kernel operator to be bounded on $L_p$ for $1<p<\infty$. It can be expressed (without details about definition) by the following theorem: 

\begin{theorem}
Let $n \in \mathbb{N}$, $p\in (1,\infty)$ and $k$ a measurable function from $\mathbb{R}^{2n}$ to $\mathbb{C}$ verifying the size and smoothness conditions. Then if the formal expression : $$Tf(x) = \int_{y\in\Rn} k(x,y)f(y)dy$$ defines a bounded operator $T$ on $L_2(\mathbb{R}^n)$, it also defines a bounded operator (still noted $T$) on $L_p(\mathbb{R}^n)$. $T$ is called a $\CZ$ operator and $k$ its kernel.
\end{theorem}

\emph{Size} and \emph{smoothness} will be defined later in a more general context. For the $p=1$ case, only weak boundedness is true in general:

\begin{theorem}\label{cascom}
With the same conditions and notations, $T$ defines an operator on $L_1(\Rn)$ and there exists a constant $C$ such that for all $f \in L_1(\Rn)$: $$ \supb{t>0}\ t\mu\{Tf > t\} \leq C\Norm{f}_1$$ where $\mu$ is the Lebesgue measure.
\end{theorem}

A motivation to show this second theorem is that it directly implies the first one by real interpolation and duality.

The ideas behind these two theorems remain valid for kernels and functions taking values in different vector spaces, which makes them a great tool to show boundedness of certain operators such as generalized Hilbert transform or Littlewood-Paley inequalities. What we prove in this paper is a generalisation of the second theorem to noncommutative integration. Furthermore, the result has already proven to be useful since it is the main ingredient used in \cite{Sukochev1} and \cite{Sukochev2} in which it is shown that  there exists a constant $c$ such that for any Lipshitz function $f$ and for any self-adjoints operators $x$ and $y$, $$\Norm{f(x) - f(y)}_{1,\infty} \leq c\Norm{f'}_{\infty}\Norm{x-y}_1.$$ This inequality does not directly express the weak $(1,1)$ boundedness of a $\CZ$ operator but is reduced to it in the mentionned papers. We also hope that our main result is a way to tackle generalisations of classical inequalities on $L_p$ whose proof relies on $\CZ$ theory. Another approach for this kind of problem is to show BMO boundedness rather than weak boundedness and conclude thanks to interpolation theory. This strategy is often easier. It first appears in \cite{JungeMeiParcet} to show the boundedness of some Fourier multipliers. And has also been applied in \cite{Xu1}, where $L_p$-boundedness of $\CZ$ operators with operator valued kernels and column valued functions is used to study Hardy spaces on quantum tori as well as in \cite{JungeParcet} where it is applied to fully noncommutative $\CZ$ operators, in quantum euclidean spaces. 

An introduction to noncommutative $L_p$-spaces can be found in \cite{PisierXu}. We will only briefly recall some basic definitions and results. A noncommutative measure space is a von Neumann algebra $\M$ equipped with a semifinite normal faithful trace $\tau$. For all $x\in\M$, we can define by the functional calculus : $$\Norm{x}_p = \tau(\md{x}^p)^{1/p}.$$ Denote $S_p = \{ x\in\M : \Norm{x}_p < \infty \}$ and $L_p(\M) = \overline{(S_p,\Norm{.}_p)}$. The elements of $L_p(\M)$ can be identified with unbounded operator affiliated to $\M$. A large part of classical integration theory still holds in this context such as Hölder's inequality, duality and interpolation.  In particular, we will need the noncommutative concept of martingales. First, a \emph{filtration} on $\M$ is an increasing sequence $(\M_n)_{n\in\Nb}$ of von Neumann subalgebras of $\M$ with weak-$\star$ dense union and such that $\tau$ restricted to each $\M_n$ remains semifinite. This guarantees the existence of conditional expectations $\E_n$ on $\M_n$ which extend to contractions from $L_p(\M)$ to $L_p(\M_n)$ for all $p \geq 1$. With this in mind, the definition of martingale is straightforward.

\subsection{Main theorem}

We will now introduce the notations that will allow us to state the main theorem of this paper. Let $(\Mt,\tau_{\Mt})$ be a noncommutative measure space and $\M$ a von Neumann subalgebra of $\Mt$ such that $\tau_{\Mt}$ restricted to $\M$ (denoted $\tau_{\M}$) is semifinite. $L_p(\M)$ is naturally included in $L_p(\Mt)$ for all $1 \leq p \leq \infty$. Let $\M'$ be the commutant of $\M$. We will consider $\N$ and $\Nt$ the von Neumann algebras of $\ast$-weakly measurable $\M$- and $\Mt$-valued functions on $\Rn$ i.e the von Neumann tensor products $\M \otimest L_{\infty}(\Rn)$ and $\Mt \otimest L_{\infty}(\Rn)$, $\tau$ will denote their natural trace (with no ambiguity since $\N$ is naturally included in $\Nt$). Let $T$ be a $\CZ$ operator associated with a kernel $k: \mathbb{R}^n \times \mathbb{R}^n \to \M' \cap \Mt$, formally given by the expression:
\begin{center}
$Tf(x) = \int_{y\in\Rn} k(x,y)f(y)dy$
\end{center}

\begin{definition} \label{def:smooth}
Say that $T$ has Lipschitz parameter $\gamma$ ($0 < \gamma \leq 1$) if for all $x$,$y$ and $z$ in $\mathbb{R}^n$ verifying $\md{x-z} \leq \frac{1}{2} \md{y-z}$, the following \emph{smoothness estimates} hold: 

\begin{align*}
\NormMt{k(x,y) - k(z,y)} &\leq \dfrac{\md{x-z}^\gamma}{\md{y-z}^{n + \gamma}},\\ 
\NormMt{k(y,x) - k(y,z)} &\leq \dfrac{\md{x-z}^\gamma}{\md{y-z}^{n + \gamma}}.
\end{align*}
\end{definition}

We will also suppose that $T$ verifies the \emph{size condition}, for all $x,y \in \Rn$:
\begin{center}
$ \NormMt{k(x,y)} \leq \dfrac{1}{\md{x-y}^n} $.
\end{center}

\begin{remark}
This size condition will never appear throughout the proof. It is however used in  remark \ref{approx:T}. Moreover, it is implied by the smoothness condition provided that $k$ goes to $0$ at infinity. 

On the contrary, the smoothness condition will be crucial to many computations. A Hörmander type condition:
\[
\supb{x\in\R^n} \Int{\md{x-y}>2\md{x-x'}}{} \md{k(x,y) - k(x',y)}dy < C
\]  
would not suffice to use the ideas of this paper. We would at least need a decay of the following type:
\[
\supb{x\in\R^n} \Int{\md{x-y}>\alpha\md{x-x'}}{} \md{k(x,y) - k(x',y)}dy < \dfrac{C}{\alpha^n},\ \alpha > 2.
\]
But we do not know if it is a strong enough condition since part of the proof of pseudo-localisation relies on pointwise estimates of the kernel.
\end{remark}

Define, for all $t > 0$ and $f \in \Nt$:
\begin{center}
$\lambda_t(f) = \tau(\{f > t\})$.
\end{center}
The main theorem we are aiming to prove in this article is the following.

\begin{theorem} \label{th:main}
There exists a constant $c_{n,\gamma}$ depending only on $n$ and $\gamma$ such that for all $\CZ$ operator $T$ with $(\M' \cap \Mt)$-valued kernel and Lipschitz parameter $0 < \gamma \leq 1$, verifying the size estimate and bounded from $L_2(\N)$ to $L_2(\Nt)$, for all $f \in L_1(\N)$:
\begin{center}
$\sup_{t>0} t\lambda_t(Tf) \leq c_{n,\gamma}\Norm{f}_1$
\end{center} 
\end{theorem}

To prove this theorem, we will need a pseudo-localisation result which constitutes the first part of this paper. This is where the most important simplifications are made compared to J. Parcet's work. In particular, our proof is elementary and does not explicitely use the size condition. We decide to directly show the result with noncommutative variables but the pseudo-localisation theorem is essentially a commutative one. The second part, which is the proof of the main theorem has been shortened and clarified thanks to more efficient organisation and computations but the underlying ideas all appear in \cite{czoperator}. In particular, we use the same decomposition which relies on Cuculescu's projections. It is a natural noncommutative counterpart to the decomposition used for the classical proof but unfortunately, new kind of terms appear which will necessitate more refined estimates and in particular, the pseudo-localisation theorem. In a third section, we show a boundedness result for singular integrals with Hilbert valued kernels and operator valued functions. It was already shown in \cite{MeiParcet} by J. Parcet and T. Mei but follows directly from theorem \ref{th:main} and Khintchine inequalities. We conclude this paper by providing an $L_p$-pseudo-localisation theorem similar to the one of P. Hytönen in \cite{Hytonen}. We would like to thank the referee for bringing the latter to our attention. The result follows mainly from our proof of $L_2$-pseudo-localisation, martingale inequalities and interpolation.

\subsection{A technical remark and notations}

The following remark allows us to manipulate the integral expression of $T$, in particular to use the smoothness condition.

\begin{remark} \label{approx:T}
For any $\CZ$ operator $T$ there exists a $\CZ$ operator $T'$ such that $Tf = T'f + gf$ for all $f$ in $L_2(\N)$ and $T'$ is the strong limit in $\B(L_2(\N))$ of operators $T_i$ given by: 
\[
T_i f(x) = \Int{y\in\R^n}{} k_i(x,y)f(y)dy,
\]
for any $f\in L_2(\N)$. Here, $k_i$ can be taken as a truncation of $k$ so that the integral makes sense and $g$ is a bounded function in $\Nt$. This fact is explained in more details in \cite{Grafakos}, proposition $8.1.11$. As a consequence, it suffices to prove the theorem for operators given by a converging integral formula on $L_2(\N)$.
\end{remark}

The proof of the main theorem will rely on the use of dyadic martingales which will require a few notations.
\begin{itemize}
\item Since we deal with cubes, the $\infty$-norm will be easier to manipulate than the euclidean one. Hence, for all $x\in\Rn$, we set $\md{x} = \Norm{x}_{\infty}$ in the remainder of the paper.  
\item $\Q$ will denote the set of all dyadic cubes and $\Q_k$ the set of dyadic cubes of edge length $2^{-k}$. Let $V_k = 2^{-nk}$ be the volume of such a cube.
\item For all $x$ in $\Rn$, $Q_{x,k}$ will denote the cube in $\Q_k$ containing $x$ and $c_{x,k}$ its center.
\item Let $(\E_k)_{k \in \mathbb{Z}}$ be the martingale associated with dyadic filtration, i.e for all $f$ :
\begin{center}
$\E_k(f)(x) = \frac{1}{V_k}\int_{Q_{x,k}}f(t)dt$
\end{center}
For convenience, we will write $f_k := \E_k(f)$ and $\Delta_k(f) := f_k - f_{k-1} =: df_k$. The filtration associated with these expectations will be denoted by $(\N_k)_{k\in\mathbb{Z}}$ where $\N_k$ is the von Neumann subalgebra of $\N$ constituted of functions that are constant on cubes of edge length $2^{-k}$.
\item For any odd positive integer $i$ and $Q$ in $\Q_k$, $iQ$ will designate the image of $Q$ by the homothety of center $c_Q$ and parameter $i$ such that $iQ$ is the union of $i^n$ cubes in $\Q_k$.
\item Notice that for all $x,y \in \mathbb{R}^{n}$ and $k\in\mathbb{Z}$, $x \in iQ_{y,k} \Leftrightarrow y \in iQ_{x,k}$. 
\end{itemize}

The notation $A \les B$ will stand for "there exists a constant $c$ depending only on $n$ and $\gamma$ such that $A \leq cB$".

In the next section, we will frequently use "polar" changes of coordinates with respect to the norm $\md{.}$ since it is more adapted to our problem. The spherical element of volume is replaced by the border of a cube which leads to a similar formula: 
\begin{center}
$\int_{\mathbb{R}^n} f(\md{x})dx = \int_{\mathbb{R}^+}2n(2r)^{n-1}f(r)dr $.
\end{center}

\subsection*{Acknowledgements}

I am very grateful to Eric Ricard for his advice and patience throughout the preparation of this article. I also thank Javier Parcet for his reading and comments on previous versions of the paper and the referee for his fruitful remarks.

\section{Pseudo-localisation}

\subsection{Theorem}

A localisation result would be of the form $\text{supp}\ Tf \approx \text{supp}\ f$ which is ideal to show weak boundedness. The theorem that follows expresses in a way the fact that singular integrals rapidly vanish outside of the support of the function to which they are applied. A simpler result of this type appears in the commutative proof in the $L_1$-context and is enough to conclude in this case. But, as mentionned before, the noncommutative case requires new tools such as the $L_2$-version of pseudo-localisation that follows. Note that the proof is written with operator-valued functions because we will need this result later but is almost a copy of the proof we had with scalar functions. So, suprisingly, the most technical part of the proof of our main theorem is purely commmutative.  

\begin{theorem} \label{th:pseudoloc}
Let $f \in L_2(\N)$ and $s \in \mathbb{N}$. For all $k\in \Z$, let $A_k$ and $B_k$ be projections in $\N_k$ such that $ A_k^{\bot}df_{k+s} = df_{k+s}B_k^{\bot} = 0$. For any odd positive integer $d$, write:
\begin{center}
$A_k = \Sum{Q\in \Q_k}{} A_Q\I_Q$, $A_Q \in \M$ and define $dA_k := \bigvee_{Q\in \Q_k} A_Q\I_{dQ}$.
\end{center}
Define $dB_k$ the same way. Let:
\[
A_{f,s} := \bigvee_{k \in \mathbb{Z}} 5A_k \text{ and } B_{f,s} := \bigvee_{k \in \mathbb{Z}} 5B_k. \addtag \label{eq:diagsize}
\]
Let $T$ be a $\CZ$ operator associated with a kernel $k$ with Lipschitz parameter $\gamma$, verifying the size condition, taking values in $\M' \cap \Mt$ and bounded from $L_2(\N)$ to $L_2(\Nt)$. Then for all $s \in \mathbb{N}$ and $f \in L_2(\N)$ we have: 
\begin{center}
$\Norm{A_{f,s}^{\bot}(Tf)}_2 \les 2^{-\frac{\gamma s}{2}}\Norm{f}_2$ and $\Norm{(Tf)B_{f,s}^{\bot}}_2 \les 2^{-\frac{\gamma s}{2}}\Norm{f}_2.$
\end{center}
\end{theorem}

Throughout the course of the proof we will often use the following lemma:

\begin{lemma}[Schur]
Let $T$ be an operator on $L_2(\N)$ given by a $\Mt$-valued kernel:
\[
Tf(x) = \int_{\mathbb{R}^n} k(x,y)f(y) dy.
\]
Let $S_1(x) = \int_{\mathbb{R}^n} \NormMt{k(x,y)} dy$ and $S_2(y) = \int_{\mathbb{R}^n} \NormMt{k(x,y)} dx$ then :
\begin{center}
$\normi{T} \leq \sqrt{\Norm{S_1}_{\infty}\Norm{S_2}_{\infty}}.$
\end{center} 
\end{lemma}

\emph{Proof}. This is not different from the commutative case (see \cite{czoperator}). Let $f \in L_2(\N)$:

\begin{align*}
\Norm{Tf}_2^2 &= \Int{\Rn}{} \Norm{\Int{\Rn}{} k(x,y)f(y)dy}_2^2 dx \\
&\leq \Int{\Rn}{} \big(\Int{\Rn}{} \Norm{k(x,y)f(y)}_2 dy \big)^2 dx \\
&\leq \Int{\Rn}{} \big(\Int{\Rn}{} \Norm{k(x,y)}_{\Mt}\Norm{f(y)}_2dy \big)^2 dx \\
&\leq \Int{\Rn}{} \Int{\Rn}{} \Norm{k(x,y)}_{\Mt} dy \Int{\Rn}{}\Norm{k(x,y)}_{\Mt}\Norm{f(y)}_2^2 dy dx \\
&\leq \Norm{S_1}_{\infty} \Int{\Rn}{} \Int{\Rn}{}\Norm{k(x,y)}_{\Mt}\Norm{f(y)}_2^2 dydx \\
&\leq \Norm{S_1}_{\infty} \Norm{S_2}_{\infty} \Int{\Rn}{} \Norm{f(y)}_2^2 dy = \Norm{S_1}_{\infty} \Norm{S_2}_{\infty} \Norm{f}_2^2
\end{align*}

\hfill{$\qed$}

It will also be important to note that by construction, for all $x,y \in \mathbb{R}^n$ such that $x \in 5Q_{k,y}$, we have $A_k(y) \leq (5A_k)(x) \leq A_{f,s}(x)$ and $A_k^{\bot}(y)df_{k+s}(y) = 0$. Consequently:
\[
A_{f,s}^{\bot}(x)df_{k+s}(y) = (5A_k)^{\bot}(x)df_{k+s}(y) = 0. \addtag \label{eq:0}
\]

We will only show the pseudo-localisation theorem in the case of left multiplications since the exact same proof can be writen for right multiplication. Another way of seeing it is that left and right multiplication are equivalent by taking adjoints. The general strategy will be to find an operator $T'$ verifying $A_{f,s}^{\bot}Tf = A_{f,s}^{\bot}T'f$ and such that we can control $\normi{T'}$ thanks to the lemma above. By remark $\ref{approx:T}$, we can also suppose that the integral defining $T$ converges. The result for any $T$ follows then directly by approximations.

\subsection{The $s$ shift}

Let $f \in L_2(\N)$, fixed throughout the proof. Let $s\in \mathbb{N}$. To make use of the construction of $A_{f,s}$ we immediately write $T = \sum_{k\in \mathbb{Z}} T \Delta_{k+s}$ where the sum converges for the strong operator topology on $\B(L_2(\N))$. We will show in this section that the constant $2^{-\gamma s}$ appears quite easily thanks to the smoothness condition when estimating the norm of $T\Delta_{k+s}f$. The tricky part will be to glue these pieces back together in the following sections. Let $k\in \Z$, note that :
$$ \int_{\Rn} k(x,c_{y,k+s-1})df_{k+s}(y)dy = \int_{\Rn}\E_{k+s-1}\big(k(x,c_{.,k+s-1})df_{k+s}(.)\big)(y)dy = 0 .$$

Therefore we can write:
\begin{align*}
A_{f,s}^{\bot}(x)Tdf_{k+s}(x) &= A_{f,s}^{\bot}(x)\int_{\Rn} k(x,y)df_{k+s}(y)dy \\
&= A_{f,s}^{\bot}(x)\int_{\Rn} (k(x,y) - k(x,c_{y,k+s-1}))df_{k+s}(y)dy \\
\intertext{$A_k^{\bot}(y)$ commutes with $k(x,y)$ and $k(x,c_{y,k+s-1})$ since $k$ is $\M'$-valued so:}
A_{f,s}^{\bot}(x)Tdf_{k+s}(x) &=  A_{f,s}^{\bot}(x)\int_{\Rn} (k(x,y) - k(x,c_{y,k+s-1}))(5A_k)^{\bot}(x)df_{k+s}(y)dy \\
\intertext{and by (\ref{eq:0}):}
&=  A_{f,s}^{\bot}(x)\int_{\Rn} \I_{x \notin 5Q_{k,y}}(k(x,y) - k(x,c_{y,k+s-1}))df_{k+s}(y)dy. 
\end{align*}
Denote by $T_k$ the operator associated with the kernel 
\[
k_k : (x,y) \mapsto \I_{x \notin 5Q_{k,y}}k(x,y) \addtag \label{def:k_k}
\]
and $T_{k,s}$ the operator associated with the kernel 
\[
k_{k,s} : (x,y) \mapsto \I_{x \notin 5Q_{k,y}}(k(x,y) - k(x,c_{y,k+s-1})). \addtag \label{def:k^s_k}
\]
It will be useful to express the result of the previous computation in terms of these operators:

\begin{lemma}
For all $k\in\Z$ and $s\in\mathbb{N}$:
\[
A_{f,s}^{\bot}T \Delta_{k+s}f = A_{f,s}^{\bot}T_k \Delta_{k+s}f = A_{f,s}^{\bot}T_{k,s} \Delta_{k+s}f \addtag \label{eq:Tk}, 
\]
and more precisely: 
\[
T_{k} \Delta_{k+s}f = T_{k,s} \Delta_{k+s}f \addtag \label{eq:Tk'}.
\]
\end{lemma}

The introduction of $T_{k,s}$ is motivated by the following lemma. 

\begin{lemma} \label{es:Tks}
For all $k\in\Z$ and $s\in\mathbb{N}$: 
\[
\normi{T_{k,s}} \les 2^{-s\gamma} \addtag \label{es:Tk'}.
\]
\end{lemma}

\emph{Proof}. The smoothness hypothesis on $T$ (see definition \ref{def:smooth}) gives:
\begin{align*}
\NormMt{k_{k,s}(x,y)} &= \I_{x \notin 5Q_{k,y}}\NormMt{k(x,y) - k(x,c_{y,k+s-1})}\\
&\leq \I_{x \notin 5Q_{k,y}}\dfrac{\md{y-c_{y,k+s-1}}^{\gamma}}{\md{y-x}^{n+\gamma}}\\
&\les \I_{x \notin 5Q_{k,y}}\dfrac{2^{-(k+s)\gamma}}{\md{y-x}^{n + \gamma}} \addtag \label{es:pointTk'}.
\end{align*}

The condition $\md{y-c_{y,k+s-1}} \leq \frac{1}{2} \md{x-y}$ is verified even for $s = 0$ as long as $x \notin 5Q_{k,y}$. This estimate is enough to apply Schur's lemma and we obtain the result by a direct computation using a "polar" change of coordinates.

\hfill{$\qed$}

\subsection{One more cancellation property}

We introduce another kernel modification whose purpose is not immediately clear but will be crucial in obtaining the estimates to apply Schur's lemma later.

\begin{property}
For all $k\in\Z$ and $s\in\mathbb{N}$ there exists an operator $S_{k,s}$ associated with a kernel $s_{k,s}$ such that the three following conditions hold:
\begin{enumerate}
\item $A_{f,s}^{\bot}T \Delta_{k+s}f = A_{f,s}^{\bot}S_{k,s} \Delta_{k+s}f$,
\item $\int_{y\in\Rn}{s_{k,s}(x,y)dy} = 0 $ for all $x\in\Rn$,
\item $\NormMt{s_{k,s}(x,y)} \les \I_{\md{x-y} > 2^{-k-1}} \dfrac{2^{-(k+s)\gamma}}{\md{y-x}^{n + \gamma}}$.
\end{enumerate}
\end{property}

\emph{Proof}. Consider $R_{k,s}$ an operator associated with the kernel $r_{k,s}$ defined by:
\begin{center}
$r_{k,s}(x,y) := \I_{\A_k}\dfrac{K(x)2^{-(k+s)\gamma}}{I_{n,\gamma}\md{y-x}^{n + \gamma}}$ 
\end{center}
where $\A_k := \{(x,y) : x\in 5Q_{y,k}, x\notin 3Q_{y,k}\}$, $K(x) = -\int_{\Rn}{k_{k,s}(x,t)dt}$ and $I_{n,\gamma} = \int_{\R^n}\I_{\A_k}(0,t){\dfrac{2^{-(k+s)\gamma}}{\md{t}^{n + \gamma}}}dt$.

We claim that $S_{k,s} := T_{k,s} + R_{k,s}$ satisfies the conditions above.

\emph{Condition $1$}. $A_{f,s}^{\bot}R_{k,s} \Delta_{k+s}f = 0$ is verified since:
\begin{center}
$\text{supp}\ r_{k,s} \subset \{(x,y) : x\in 5Q_{y,k}\}$ i.e $r_{k,s} = \I_{x\in 5Q_{y,k}}r_{k,s}$.
\end{center}
Indeed, the shift section shows that to compute $A_{f,s}^{\bot}R_{k,s} \Delta_{k+s}f$, and in particular $(5A_k)^{\bot}R_{k,s} \Delta_{k+s}f$, we might as well replace the kernel $r_{k,s}$ by $\I_{x \notin 5Q_{k,y}}r_{k,s} = \I_{x \notin 5Q_{k,y}}\I_{x\in 5Q_{y,k}}r_{k,s} = 0$.

\emph{Condition $2$}. It is direct by definition of $K$ and $I_{n,\gamma}$. Let $x\in\Rn$:
\begin{center}
$\int_{y\in\Rn} s_{k,s}(x,y)dy  = \int_{\Rn}{k_{k,s}(x,y)dy} + K(x) = 0$.
\end{center}

\emph{Condition $3$}. It suffices to show that $\NormMt{K(x)}$ is bounded by a constant depending only on $n$ and $\gamma$. By a "polar" change of coordinates, there exists a constant $c_n$ such that:
\begin{center}
$\NormMt{K(x)} \leq \int_{\md{t} >4.2^{-k}}{\dfrac{2^{-(k+s)\gamma}}{\md{t}^{n + \gamma}}}dt = c_n2^{-\gamma s}$.
\end{center}
The bound does not depend on $x$, as we needed.

\hfill{$\qed$}

\subsection{The decomposition}

We are ready to decompose $T$ into operators whose norms in $\B(L_2(\N))$ are controlled. Fix $s\in\mathbb{N}$, write: 

\begin{align*}
A_{f,s}^{\bot}Tf &= A_{f,s}^{\bot}\sum_{k\in \mathbb{Z}} T\Delta_{k+s}f = A_{f,s}^{\bot}\sum_{k \in \mathbb{Z}} (T_k + R_{k,s})\Delta_{k+s}f \\
&= A_{f,s}^{\bot}\sum_{k\in \mathbb{Z}}\sum_{i\in \mathbb{Z}} \Delta_{k+i}(T_k + R_{k,s})\Delta_{k+s}f \\
\end{align*} 
Note that for $i \geq 1$, $(5A_k)^{\bot}$ commutes with $\Delta_{k+i}$ and recall that the first point of $R_{k,s}$'s construction implies that $(5A_k)^{\bot}R_{k,s}\Delta_{k+s}f = 0$, then:
\begin{align*}
A_{f,s}^{\bot}\Delta_{k+i}(T_k + R_{k,s})\Delta_{k+s}f &= A_{f,s}^{\bot}\Delta_{k+i}(5A_k)^{\bot}(T_k + R_{k,s})\Delta_{k+s}f \\
&= A_{f,s}^{\bot}\Delta_{k+i}T_k\Delta_{k+s}f
\end{align*}
For any $i \in \mathbb{Z}$, we have:
\begin{align*}
A_{f,s}^{\bot}\Delta_{k+i}(T_k + R_{k,s})\Delta_{k+s}f &= A_{f,s}^{\bot}\Delta_{k+i}(T_{k,s} + R_{k,s})\Delta_{k+s}f \\
&= A_{f,s}^{\bot}\Delta_{k+i}S_{k,s}\Delta_{k+s}f
\end{align*}
Now we can write our final decomposition:
\[
A_{f,s}^{\bot}Tf = A_{f,s}^{\bot}\big(\Sum{i=0}{\infty} \Phi_i f + \Sum{i=1}{\infty}\Psi_i f\big) \addtag{} \label{Decomposition}
\]
where: 
\begin{center}
$\Phi_i = \Sum{k\in \Z}{} \Delta_{k+i}T_k\Delta_{k+s}$ and
$\Psi_i = \Sum{k\in \Z}{} \Delta_{k-i}S_{k,s}\Delta_{k+s}$
\end{center}

We have then reduced our pseudo-localisation theorem to the following proposition.

\begin{property}\label{es:phi_i psi_i}
The following estimates hold for $i\geq 0$:
\begin{center}
$\normi{\Phi_i} \les 2^{-\gamma \frac{i+s}{2}}$  and $\normi{\Psi_i} \les \sqrt{1+i}2^{-\gamma \frac{i+s}{2}}$
\end{center}
\end{property}

The proof of this result will occupy the next two parts. 

\subsection{Estimate for $\Phi_i$}

Let $i \geq 0$. Note that $T^*$ is also a $\CZ$ operator, associated with the kernel $k^{\star}$ where: $k^{\star}(x,y) = k(y,x)^*$.

So $k^{\star}$ satisfies the smoothness condition of parameter $\gamma$, which means by $(\ref{es:Tk'})$ that:
\begin{center}
$\normi{(T^*)_{k,s}} \les 2^{-s\gamma}$.
\end{center}
 
Consequently :

\begin{center}
$\normi{\Delta_{k+i}T_k} = \normi{T^*_k\Delta_{k+i}} = \normi{(T^*)_{k,i} \Delta_{k+i}} \les 2^{-\gamma i}$,
\end{center}
where we used that $(T_k)^* = (T^*)_k$ and $(\ref{eq:Tk'})$. 

By orthogonality of the $\Delta_k$, we have on one hand:
\begin{center}
$\normi{\Phi_i} \leq \sup_k \normi{\Delta_{k+i}T_k\Delta_{k+s}} \leq \sup_k \normi{\Delta_{k+i}T_k} \les 2^{-\gamma i}$,
\end{center}
and on the other hand :
\begin{center}
$\normi{\Phi_i} \leq \sup_k \normi{\Delta_{k+i}T_k\Delta_{k+s}} \leq \sup_k \normi{T_k\Delta_{k+s}} \les 2^{-\gamma s}$.
\end{center}
By combining the two, $\normi{\Phi_i} \les 2^{-\gamma\frac{s+i}{2}}$.

\subsection{Estimate for $\Psi_i$}

\begin{lemma} \label{es:laid}
For all $x \in\Rn$, $k\in\Z$ and $i\in\mathbb{N}$, the following estimate holds:
\begin{center}
$\int_{t\in Q_{x,k-i}}\int_{y\in Q_{x,k-i}^c} \I_{\md{t-y}>2^{1-k}} \dfrac{1}{\md{t-y}^{n + \gamma}}dydt \les (1+i)2^{(\gamma - n)(k-i)} $.
\end{center}
\end{lemma}

\emph{Proof}. Fix $x \in\Rn$, $k\in\Z$ and $i\in\mathbb{N}$. Denote $Q = Q_{x,k-i}$, $c$ the center of $Q$ and:
\begin{center}
$X = \int_{t\in Q}\int_{y\in Q^c} \I_{\md{t-y}>2^{1-k}} \dfrac{1}{\md{t-y}^{n + \gamma}}dydt$.
\end{center} 

For every $t$ in $Q$ let $\delta_t$ be the distance from $t$ to $Q^c$ and notice that for any $t$: 

\begin{center}
$\int_{Q^c}\I_{\md{t-y}>2^{1-k}}\dfrac{1}{\md{t-y}^{n+\gamma}}dy \leq \int_{\md{y-t}>\delta_t}\I_{\md{t-y}>2^{1-k}}\dfrac{1}{\md{t-y}^{n+\gamma}}dy =: f(\delta_t)$
\end{center}
Indeed, the term on the right only depends on $\delta_t$. For $r \leq 2^{-(k-i+1)}$:
\begin{center}
$\{ t \in Q : \delta_t = r\} = \{t \in Q : \md{t - c} = 2^{-(k-i+1)} - r =: r'\}$,
\end{center}
So by a "polar" change of coordinates:
\begin{center}
$X \leq \int_{0}^{2^{-(k-i+1)}} 2n(2r')^{n-1}f(r)dr \les 2^{-(k-i)(n-1)}\int_{0}^{2^{-(k-i+1)}} f(r)dr$.
\end{center}
By a direct computation:
\[
f(r) \les \left\{
    \begin{array}{ll}
        r^{-\gamma} & \mbox{for all}\ r \\
        2^{\gamma(k-1)} & \mbox{if moreover}\ r \leq 2^{-(k-1)}
    \end{array}
\right.
\]
For $0 < \gamma < 1$ we only need the first estimate:

\begin{align*}
X &\les 2^{-(k-i)(n-1)}\int_{0}^{2^{-(k-i+1)}} f(r)dr &\les&\ 2^{-(k-i)(n-1)}\int_{0}^{2^{-(k-i+1)}} r^{-\gamma}dr \\
&\les 2^{-(k-i)(n-1)} 2^{-(k-i+1)(1-\gamma)} &\les&\ 2^{-n(k-i)}2^{\gamma(k-i)} 
\end{align*}

For $\gamma = 1$ we have to decompose the integral according to the distinction made above when estimating $f$:
\begin{align*}
\int_{0}^{2^{-(k-i+1)}} f(r)dr =& \int_{0}^{2^{-(k-1)}} f(r)dr + 
\int_{2^{-(k-1)}}^{2^{-(k-i+1)}} f(r)dr \\
&\les 2^{-(k-1)}2^{\gamma(k-1)} + \int_{2^{-(k-1)}}^{2^{-(k-i+1)}} \dfrac{1}{r} dr \\
&\les 1 + i 
\end{align*}
Therefore:
\begin{center}
$X \les 2^{-(k-i)(n-1)}(1 + i) = (1+i)2^{-n(k-i)}2^{\gamma(k-i)}$.
\end{center}

\hfill{$\qed$}

\begin{property}
For all $k$ in $\mathbb{Z}$ and $i > 0$, we have the following estimate:
\begin{center}
$\normi{\mathcal{E}_{k-i}S_{k,s}} \les \sqrt{1+i}2^{-\gamma(s+i/2)}$
\end{center}
\end{property}

\emph{Proof.} For any function $g : \mathcal{E}_{k-i}S_{k,s}g(x) = \int_{Q_{x,k-i}}\frac{1}{V_{k-i}}\int_{\Rn}s_{k,s}(t,y)g(y)dydt$. So $\mathcal{E}_{k-i}S_{k,s}$ corresponds to the kernel:
\begin{center}
$E : (x,y) \mapsto \frac{1}{V_{k-i}}\int_{Q_{x,k-i}}s_{k,s}(t,y)dt = \frac{1}{V_{k-i}}\int_{Q_{x,k-i}^c}s_{k,s}(t,y)dt $,
\end{center} 
where we used the cancellation property on $S_{k,s}$. We are looking to estimate both integrals in order to apply Schur's lemmma. Fix $x$:

\begin{align*}
\int_{\Rn}\NormMt{E(x,y)}dy =& \int_{Q_{x,k-i}}\NormMt{E(x,y)}dy + \int_{Q_{x,k-i}^c}\NormMt{E(x,y)}dy \\
=& \int_{Q_{x,k-i}}\NormMt{\int_{Q_{x,k-i}^c}\frac{1}{V_{k-i}}s_{k,s}(t,y)dt}dy\\ +& \int_{Q_{x,k-i}^c}\NormMt{\int_{Q_{x,k-i}}\frac{1}{V_{k-i}}s_{k,s}(t,y)dt}dy \\
\les& \frac{1}{V_{k-i}} \int_{Q_{x,k-i}}\int_{Q_{x,k-i}^c}\I_{\md{t-y}>2^{1-k}}\dfrac{2^{-\gamma(k+s)}}{\md{t-y}^{n+\gamma}}dtdy \\
\intertext{Using lemma $\ref{es:laid}$ :}
\int_{\Rn}\NormMt{E(x,y)}dy \les&\ \frac{1}{V_{k-i}}(1+i)2^{-n(k-i)}2^{\gamma(k-i)}2^{-\gamma(k+s)} \\
\les&\ (1+i)2^{-\gamma(i+s)} 
\end{align*}

The other estimate, with $y$ fixed, is straightforward:

\begin{align*}
\int_{\Rn}\NormMt{E(x,y)}dx =& \int_{\Rn}\NormMt{\dfrac{1}{V_{k-i}}\int_{Q_{x,k-i}}s_{k,s}(t,y)dt}dx \\
\leq& \int_{\Rn}\dfrac{1}{V_{k-i}}\int_{Q_{x,k-i}}\NormMt{s_{k,s}(t,y)}dtdx \\
=& \int_{\Rn}\NormMt{s_{k,s}(t,y)}dt \\
\les&\ 2^{-\gamma s}
\end{align*}

The proposition follows directly from the two previous computations and Schur's lemma. 

\hfill{$\qed$}

We now have:
\begin{center}
$\normi{\Delta_{k-i}S_{k,s}\Delta_{k+s}} \leq \normi{\mathcal{E}_{k-i}S_{k,s}} \les \sqrt{1+i}\ 2^{-\gamma\frac{s+i}{2}}$. 
\end{center}
We conclude once more thanks to the orthogonality between the $\Delta_{k-i}S_{k,s}\Delta_{k+s}$: 
\begin{center}
$\normi{\Psi_i} \les \sup_k \normi{\Delta_{k-i}S_{k,s}\Delta_{k+s}} \les \sqrt{1+i}\ 2^{-\gamma \frac{s+i}{2}}$
\end{center}

This concludes the proof of the pseudo-localisation theorem.

\section{Proof of the main theorem}

\subsection{The good and the bad functions}

The idea of the following decomposition comes from the classical proof of the weak type inequality for singular integrals. It has been noticed that the commutative decomposition can be expressed in terms of martingales which is well suited for a translation in the noncommutative setting. However, even with this idea, the construction is not immediate and the estimates are more difficult to obtain due to the appearance of new "off-diagonal" terms.

Fix $t > 0$ and $f \in L_1(\N)$. We will suppose that $f$ is positive to avoid unnecessary computations. This is possible because $f$ can always be written as: $$f = f_1 - f_2 + if_3 - if_4$$ with $\Norm{f_k}_1 \leq \Norm{f}_1$ for $k = 1,2,3,4$. Suppose also that $\{x \in \Rn : f(x) \neq 0 \}$ is bounded (in other words, $f$, considered as a function, has a compact support), and that $f$ is in $L_2(\N) \cap \N$. We can make these assumptions since the functions satisfying them form a dense subspace $\mathcal{S}$ of $L_1(\N)$. Once the theorem is proven, $T$ can be defined on $L_1(\N)$ as the only bounded extension of its restriction to $\mathcal{S}$.  

Denote by $(f_n)_{n \in \mathbb{Z}}$ the martingale associated with $f$ and the filtration $(\N_n)_{n\in\mathbb{Z}}$. The main tools to decompose $f$ into a good and a bad part are Cuculescu's projections:

\begin{theorem}[Cuculescu]
Let $x = (x_n)$ be a bounded positive $L_1$-martingale and $t \geq 0$. Then there exists a decreasing sequence $(q_n)$ of projections in $\N$ such that for every $n \geq 1$:
\begin{enumerate}
\item $q_n \in \N_n$.
\item $q_n$ commutes with $q_{n-1}x_nq_{n-1}$.
\item $q_nx_nq_n \leq t$.
\item moreover, if $q = \bigwedge q_n$ then : 
\begin{center}
$qx_nq \leq t$ for $n \geq 1$ and $\tau(q^{\bot}) \leq \dfrac{\Norm{x}_1}{t}$.
\end{center}
\end{enumerate}
\end{theorem}

The boundedness hypothesis on $f$ and its support imply that there exists $n_0$ such that for all $n \leq n_0$, $f_n \leq t$. Without loss of generality, we will suppose that $n_0 = 0$. From now on, let $q_n$ denote the projections given by Cuculescu's theorem, associated with $t$ and $(f_n)_{n\geq 0}$. To complete this definition let $q_n = 1$ for all $n \leq 0$. Notice that for all $n \in \Z$, $q_nf_nq_n \leq t$.

Define :
\begin{center}
$\forall n\in \Z, p_n = q_{n-1} - q_n$ and $p_{\infty} = q$.
\end{center}
Let $\Zb = \Z \cup \{\infty\}$. By definition, $p_n \in \N_n$ and:
\begin{center}
$\Sum{n\in \Zb}{}p_n = 1$.
\end{center}
Which allows us to define the good and bad parts as follows :
\begin{center}
$g = \Sum{i\in \Zb}{}\Sum{j\in \Zb}{} p_if_{i\vee j}p_j$ and $b = \Sum{i\in \Zb}{}\Sum{j\in \Zb}{} p_i(f-f_{i\vee j})p_j$.
\end{center}

By properties of the distribution function (see \cite{FackKosaki}):
\begin{center}
$\lambda_t(Tf) \les \lambda_{t/2}(Tg) + \lambda_{t/2}(Tb)$
\end{center}
So it suffices to prove estimates of the form: $$\lambda_t(Tx) \les \dfrac{\Norm{f}_1}{t}$$
for both $x = g$ and $x = b$. This is the purpose of the next sections.

\begin{remark}
The same formula for $b$ and $g$ works in the commutative case except that only the diagonal terms are non zero. This explains why the commutative proof can only be repeated for the diagonal terms and "shifted" diagonals. What makes this decomposition work is that, due to the pseudo-localisation lemma, the estimates for "shifted" diagonals get exponentially better as the shift increases. 
\end{remark}

We will use the following notation: for all $k$ and $Q \in \Q_k$, $p_Q := p_k(x)$ for any $x \in Q$. We will also need two lemmas which are directly deduced from the construction.

\begin{lemma} \label{ineq:fk}
For all $k \in \Z$, we have : $f_{k+1} \leq 2^{n}f_k$.

\end{lemma}

\emph{Proof}. This is straightforward from the definition of $f_k$ and positivity of $f$. Let $x\in \Rn$:

\begin{center}
$2^n f_k(x) = \dfrac{2^n}{V_k} \int_{Q_{x,k}} f(t)dt \geq \dfrac{2^n}{V_k} \int_{Q_{x,k+1}} f(t)dt = f_{k+1}(x)$.
\end{center}

\hfill{$\qed$}

\begin{lemma} \label{zeta}
Let $d$ be an odd positive integer. Define:
\begin{center}
$\zeta_d = \big(\bigvee_{Q\in \Q} p_Q\I_{dQ}\big)^{\bot}$.
\end{center}
Then :
\begin{enumerate}
\item $\tau(\zeta_d^{\bot}) \leq d^n\dfrac{\Norm{f}_1}{t}$
\item For all cubes $Q \in \Q$, we have the following cancellation property:
\begin{center}
$x \in dQ \Rightarrow \zeta_d(x)p_Q = p_Q\zeta_d(x) = 0$.
\end{center}
\end{enumerate}
\end{lemma}

\emph{Proof}. The first estimate is a consequence of Cuculescu's inequality:
\begin{center}
$\tau(\zeta_d^{\bot}) \leq \Sum{Q\in \Q}{} d^n\tau(p_Q\I_Q) = d^n \Sum{k=1}{\infty}\tau(p_k) = d^n \tau(q^{\bot}) \leq d^n\dfrac{\Norm{f}_1}{t}$.
\end{center}

Let $Q \in \Q$ and $x \in dQ$. By construction, $\zeta_d^{\bot}(x) \geq p_Q$ so $\zeta_d(x) \leq p_Q^{\bot}$. This concludes the proof.

\hfill{$\qed$}

From now on, we fix $d = 5$ and denote $\zeta_d$ by $\zeta$.

\begin{remark}
This projection $\zeta$ is to be thought as a dilatation of the support of the bad function which already plays a crucial role in the commutative setting.
\end{remark}

\subsection{Estimate for the bad function}

The strategy of proof is to write:
\begin{center}
$Tb = \zeta Tb \zeta + (1 - \zeta) Tb \zeta + \zeta Tb (1 - \zeta) + (1 - \zeta) Tb (1 - \zeta)$.
\end{center}
Therefore, lemma $\ref{zeta}$ and Tchebychev's inequality give:
\begin{center}
$\lambda_t(Tb) \les \tau(1 - \zeta) + \lambda_t(\zeta Tb \zeta) \les \dfrac{\Norm{f}_1}{t} + \dfrac{\Norm{\zeta Tb \zeta}_1}{t}$.
\end{center}

The estimate for the bad function is now reduced to the following proposition:
\begin{property} \label{es:bad}
We have the estimate: $\Norm{\zeta Tb \zeta}_1 \les \Norm{f}_1$.
\end{property}

The proof will require three intermediate lemmas. 

Define, for all $i,j \in \Z$ : $b_{i,j} = p_i(f-f_{i\vee j})p_j$.

\begin{lemma} \label{es:l1:b}
For all $s\in\Z$ :$\Sum{i-j = s}{} \Norm{b_{i,j}}_1 \les \Norm{f}_1$
\end{lemma}

\emph{Proof}. Let $i,j\in\Z$ :
\begin{align*}
\Norm{b_{i,j}}_1 &\leq \Norm{p_ifp_j}_1 + \Norm{p_if_{i\vee j}p_j}_1 \\
\intertext{By Holder's inequality :}
\Norm{b_{i,j}}_1 &\leq \Norm{f^{1/2}p_i}_2\Norm{f^{1/2}p_j}_2 + \Norm{f_{i\vee j}^{1/2}p_i}_2\Norm{f_{i\vee j}^{1/2}p_j}_2 \\
&\leq \dfrac{1}{2}\big(\Norm{p_ifp_i}_1 + \Norm{p_jfp_j}_1 + \Norm{p_if_{i\vee j}p_i}_1 + \Norm{p_jf_{i\vee j}p_j}_1 \big)\\
&\leq \tau(p_if) + \tau(p_jf)
\end{align*}

Consequently, for all $s\in\Z$:

\begin{center}
$\Sum{i-j = s}{} \Norm{b_{i,j}}_1 \leq 2\Sum{i\in\Z}{} \tau(p_if) \les \Norm{f}_1$
\end{center}

\hfill{$\qed$}

\begin{lemma} \label{eq:cancel:bij}
The following cancellation properties hold: 
\begin{itemize}
\item for all $i,j\in\mathbb{Z}$ and $Q\in \Q_{i\vee j}$: $\int_Q b_{i,j} = 0$;
\item for all $x,y\in\Rn$ such that $y \in 5Q_{x,i\wedge j}$: $\zeta(x)b_{i,j}(y)\zeta(x) = 0$.
\end{itemize}
\end{lemma}

\emph{Proof}. The first point is straightforward. Since $Q \in \Q_{i\vee j}$:
$\int_Q b_{i,j} = \int_Q \E_{i\vee j}(b_{i,j}) = 0$.

The second one is a consequence of the construction of $\zeta$ expressed in property $\ref{zeta}$. Indeed, for all $x,y\in\Rn$ and $k\in\mathbb{Z}$ such that $y \in 5Q_{x,k}$ we know that $\zeta(x)p_k(y) = p_k(y)\zeta(x) = 0$. Recall that $b_{i,j} = p_i(f-f_{i\vee j})p_j$, it is now clear that $\zeta(x)b_{i,j}(y)\zeta(x) = 0$ for $y \in 5Q_{x,k}$ and $k = i,j$ which concludes the proof of the lemma. Note that $Q_{x,k} \subset Q_{x,k'}$ for $k \geq k'$, so we do not lose anything by taking $k = i\wedge j$.

\hfill{$\qed$}

The following lemma is the core of the bad function estimate, it relies on a computation which allows us to make use of the smoothness condition.

\begin{lemma} \label{es:zetaTbij}
For all $i,j \in \Z$ : $\Norm{\zeta Tb_{i,j} \zeta}_1 \les 2^{-\md{i-j}\gamma}\Norm{b_{i,j}}_1$.

\end{lemma}

\emph{Proof}. Note that $b_{i,j}$ is in $L_2(\N)$. Fix $i,j\in \Z$ and $x \in \Rn$, recall that $k(x,y)$ commutes with $\zeta(x)$ since $\zeta$ takes values in $\M$ and $k$ in $\M'$:

\begin{align*}
\zeta(x)Tb_{i,j}(x)\zeta(x) &= \int_{y\in \Rn} k(x,y)\zeta(x)b_{i,j}(y)\zeta(x) dy \\
&= \int_{y \in 5Q_{x,i\wedge j}^c} \big(k(x,y) - k(x,c_{y,i\vee j})\big) \zeta(x)b_{i,j}(y)\zeta(x) dy
\end{align*}
where we used both cancellation properties of $b_{i,j}$ (lemma $\ref{eq:cancel:bij}$) the first one to make the term $k(x,c_{y,i\vee j})$ appear and the second one to reduce the domain of integration. 
Therefore, the smoothness condition (definition \ref{def:smooth}) applies and gives :

\begin{align*}
\Norm{\zeta(x)Tb_{i,j}(x)\zeta(x)}_1 &\leq \int_{y \in 5Q_{x,i\wedge j}^c} \NormMt{k(x,y) - k(x,c_{y,i\vee j})} \Norm{b_{i,j}(y)}_1 dy\\
&\leq \int_{y \in \Rn} \I_{y\notin 5Q_{x,i\wedge j}} \dfrac{2^{-\gamma(i\vee j)}}{\md{x-y}^{n+\gamma}} \Norm{b_{i,j}(y)}_1 dy \\
\end{align*}
It follows that :
\begin{align*}
\Norm{\zeta Tb_{i,j} \zeta}_1 &= \int_{x\in \Rn} \Norm{\zeta(x)Tb_{i,j}(x)\zeta(x)}_1 dx \\
&\leq \int_{y\in \Rn} \int_{x \in \Rn} \I_{x\notin 5Q_{y,i\wedge j}}\dfrac{2^{-\gamma(i\vee j)}}{\md{x-y}^{n+\gamma}} \Norm{b_{i,j}(y)}_1 dxdy \\
&\les \int_{y\in \Rn} 2^{-\gamma (i\vee j - i\wedge j) } \Norm{b_{i,j}(y)}_1 dy \\
&\les 2^{-\gamma \md{i-j}} \Norm{b_{i,j}}_1
\end{align*}

\hfill{$\qed$}

\emph{Proof of property $\ref{es:bad}$}. The only thing left to do is to glue the pieces together thanks to lemmas $\ref{es:l1:b}$ and $\ref{es:zetaTbij}$:
\begin{align*}
\Norm{\zeta Tb \zeta}_1 &\leq \Sum{i,j \in\Z}{} \Norm{\zeta Tb_{i,j}\zeta}_1 
\leq \Sum{s\in \Z}{} \Sum{i-j=s}{} \Norm{\zeta Tb_{i,j}\zeta}_1 \\
&\leq \Sum{s\in \Z}{} 2^{-\gamma\md{s}}\Sum{i-j=s}{} \Norm{b_{i,j}}_1 
\les \Sum{s\in \Z}{} 2^{-\gamma\md{s}} \Norm{f}_1 \les \Norm{f}_1. 
\end{align*}

\subsection{Estimate for the good function}

This one is more involved and requires the $L_2$-pseudo-localisation theorem. The same trick as for the bad function allows us to write:

\begin{center}
$\lambda_t(Tg) \les \tau(1 - \zeta) + \lambda_t(\zeta Tg \zeta) \les \dfrac{\Norm{f}_1}{t} + \dfrac{\Norm{\zeta Tg \zeta}_2^2}{t^2}$.
\end{center}

Define the diagonal, left and right parts of $g$ as follows:
\begin{center}
$g^{(d)} = \Sum{i\in\Zb}{} p_if_ip_i$ , $g^{(l)} = \big(\Sum{i<j\in\Z}{} p_if_jp_j\big) + q^{\bot}fq$ and $g^{(r)} = \big(\Sum{i<j\in\Z}{} p_jf_jp_i\big) + qfq^{\bot}$.
\end{center}

Note that the estimate for the good function can easily be deduced from the following property:

\begin{property} \label{es:good}
The following estimates hold :
\begin{center}
$\Norm{Tg^{(d)}}_2^2 \les t\Norm{f}_1$ , $\Norm{\zeta Tg^{(l)}}_2^2 \les t\Norm{f}_1$ and $\Norm{Tg^{(r)} \zeta}_2^2 \les t\Norm{f}_1$
\end{center}
\end{property}

\emph{Proof of the $g^{(d)}$ estimate}. Since $T$ is bounded on $L_2(\N)$ it sufices to prove that $\Norm{g^{(d)}}_2^2 \les t\Norm{f}_1$.

Notice that $g^{(d)}$ is positive.
\begin{center}
$\Norm{g^{(d)}}_1 = \Sum{i\in \Zb}{} \tau(p_if_ip_i) = \Sum{i\in \Zb}{} \tau(p_if) = \tau(\Sum{i\in \Zb}{} p_if) = \Norm{f}_1$.
\end{center} 

By orthogonality, we have $\Norm{g^{(d)}}_{\infty} = \supb{k\in \Zb}\Norm{p_kf_kp_k} \leq 2^nt $. Indeed, for $k < \infty$, we have, by Lemma $\ref{ineq:fk}$:
\begin{center}
$p_kf_kp_k \leq 2^n p_kf_{k-1}p_k \leq 2^n q_{k-1}f_{k-1}q_{k-1} \leq 2^n t$
\end{center}

For $k= \infty$, reasoning in $L_2(\N)$: $t - qfq  = \Lim{k \to \infty} t - qf_kq \geq 0$ since $L_2(\N)^+$ is closed. So $\Norm{g_d}_{\infty} \leq 2^n t$.

We conclude by Hölder's inequality:
\begin{center}
$\Norm{g_d}_2^2 \leq \Norm{g_d}_1\Norm{g_d}_{\infty} \leq 2^nt\Norm{f}_1$.
\end{center}

\hfill{$\qed$}

\emph{Proof of the $g^{(l)}$ estimate}. This will conclude the proof of proposition $\ref{es:good}$ since the argument for $g^{(r)}$ is similar.

\begin{lemma}
We have the following expression for $g^{(l)}$:
\begin{center}
$g^{(l)} = \Sum{k = 1}{\infty} \Sum{s=1}{\infty} p_kdf_{k+s}q_{k+s-1} =: \Sum{k=1}{\infty} \Sum{s=1}{\infty} g_{s,k}^{(l)}$.
\end{center}
\end{lemma}

\emph{Proof}. This is obtained by an Abel's transform. First, fix $j_0$:
\begin{align*}
\Sum{i<j\leq j_0}{} p_if_jp_j &= \Sum{i<j\leq j_0}{} p_if_j(q_{j-1}-q_j) = \Sum{i\leq j<j_0}{}p_idf_{j+1}q_j - \Sum{i<j_0}{}p_if_{j_0}q_{j_0} \\
&= \Sum{i \leq j<j_0}{}p_idf_{j+1}q_j - q_{j_0-1}^{\bot}f_{j_0}q_{j_0}
\end{align*}

Letting $j_0$ go to infinity, we obtain:
\begin{center}
$\Sum{i<j \in \Z}{} p_if_jp_j = \big( \Sum{i \leq j \in \Z}{}p_idf_{j+1}q_j \big) - q^{\bot}fq$,
\end{center}
which is exactly what we needed.

\hfill{$\qed$}

\begin{property}\label{es:gs}
Define, for all $s \geq 1$:
\begin{center}
$g_s^{(l)} = \Sum{k=1}{\infty} p_kdf_{k+s}q_{k+s-1}$
\end{center}
Then: $\Norm{\zeta Tg_s^{(l)}}_2^2 \les 2^{-\gamma s} t\Norm{f}_1$.
\end{property}

This is, as for the bad function, the core of the proof and will require some work.

\begin{lemma}
We have the following properties:
\begin{enumerate}
\item $\Delta_{k+s}(g_s^{(l)}) = g_{s,k}^{(l)}$.
\item Fix $s$, the $g_{s,k}^{(l)}$ are orthogonal in $L_2(\N)$.
\item $\Norm{g_{s}^{(l)}}_2^2 \les t \Norm{f}_1$.
\end{enumerate}
\end{lemma}

\emph{Proof}. $1$ is straightforward and $2$ follows directly from $1$ since martingale differences are always orthogonal. Let us prove $3$.
\begin{align*}
\Norm{g_{s,k}^{(l)}}_2^2 &\leq 2(\Norm{p_kf_{k+s}q_{k+s-1}}_2^2 + \Norm{p_kf_{k+s-1}q_{k+s-1}}_2^2) \\
&\leq 2\tau(p_kf_{k+s}q_{k+s-1}f_{k+s}p_k) + 2\tau(p_kf_{k+s-1}q_{k+s-1}f_{k+s-1}p_k)
\end{align*}

By Cuculescu's theorem and recalling that $f_i \leq 2^nf_{i-1}$ from lemma $\ref{ineq:fk}$:
\begin{align*}
\Norm{f_{k+s}^{1/2}q_{k+s-1}f_{k+s}^{1/2}}_{\infty} &= \Norm{q_{k+s-1}f_{k+s}q_{k+s-1}}_{\infty} \les t \\
\Norm{f_{k+s}^{1/2}q_{k+s-1}f_{k+s-1}^{1/2}}_{\infty} &= \Norm{q_{k+s-1}f_{k+s-1}q_{k+s-1}}_{\infty} \leq t
\end{align*}

Therefore, for all $s > 0$:
\begin{center}
$\Norm{g_s^{(l)}}_2^2 = \Sum{k>0}{} \Norm{g_{k,s}^{(l)}}_2^2 \les \Sum{k>0}{} t\tau(p_kf) = t\Norm{f}_1$.
\end{center}
\hfill{$\qed$}

\begin{lemma}
We have the following estimate:
\begin{center}
$\Norm{\zeta Tg_s^{(l)}}_2 \les 2^{-\gamma s/2} \Norm{g_s^{(l)}}_2$
\end{center}
\end{lemma}

\emph{Proof}. This is where we apply pseudo-localisation i.e theorem $\ref{th:pseudoloc}$ to $g_s^{(l)}$. Using the notations introduced for this theorem, we can take $A_k = p_k$ since $p_k^{\bot}g_{s,k}^{(l)} = 0$. Then:
\begin{center}
$A_{f,s} = \bigvee_{k>0}5p_k = \zeta^{\bot}$.
\end{center}
The theorem gives:
\begin{center}
$\Norm{A_{f,s}^{\bot} Tg_s^{(l)}}_2 \les 2^{-\gamma s/2}\Norm{g_s^{(l)}}_2$
\end{center}
which is exactly the expected estimate.

\hfill{$\qed$}

The proposition $\ref{es:gs}$ is clear from the two previous lemmas.

It follows that :
\begin{align*}
\Norm{\zeta Tg^{(l)}}_2^2 &\leq \big( \Sum{s=1}{\infty} \Norm{\zeta Tg_s^{(l)}}_2 \big)^2 
\leq \big( \Sum{s=1}{\infty} 2^{-\frac{\gamma s}{2}}\Norm{g_s^{(l)}}_2 \big)^2 \\
&\les \big( \Sum{s=1}{\infty} 2^{-\frac{\gamma s}{2}}\sqrt{t\Norm{f}_1} \big)^2 
\les t\Norm{f}_1
\end{align*}

which is the expected estimate for $g^{(l)}$ and concludes the proof of the main theorem.

\section{Application to Hilbert-valued kernels}

Let $(\M,\tau)$ be a noncommutative measure space, $\N = L_{\infty}(\mathbb{R}) \overline{\otimes} \M$. Let $T$ be a $\CZ$ operator associated to a kernel $k$ taking values in $\ell^2$. For all $1\leq p \leq \infty$, we have a multiplication from $\ell^2 \times L_p(\M)$ to $L_p(\M)^{\mathbb{N}}$: for all $h = (h_i)_{i\in \mathbb{N}} \in \ell^2$ and $x \in \M$, define $hx = (h_ix)_{i\in \mathbb{N}}$. 

There are different natural norms on $L_p(\M)^{\mathbb{N}}$ that make this multiplication continuous. We will be interested in (see \cite{PisierXu}):
\begin{center}
$\Norm{x}_{C_p(\M)} = \Norm{\big( \Sum{k=0}{\infty} x_k^*x_k \big)^{1/2}}_p$ , $\Norm{x}_{R_p(\M)} = \Norm{\big( \Sum{k=0}{\infty} x_kx_k^* \big)^{1/2}}_p$ 
\end{center} 
and 
$$
\Norm{x}_{RC_p(\M)} = \left\{
    \begin{array}{ll}
        \inf\{\Norm{y}_{C_p(\M)} + \Norm{z}_{R_p(\M)} : y + z = x\} & \mbox{if } p < 2 \\
        \max(\Norm{x}_{C_p(\M)},\Norm{x}_{R_p(\M)}) & \mbox{if } p \geq 2.
    \end{array}
\right.
$$
Denote by $C_p(\M)$, $R_p(\M)$ and $RC_p(\M)$ the associated subspaces of $L_p(\M)^{\mathbb{N}}$. Likewise, define $C_p(\N)$, $R_p(\N)$ and $RC_p(\N)$. Note that $C_2(\M) = R_2(\M) = RC_2(\M) = \ell^2(L_2(\M))$. 

Thanks to the multiplication we defined, $T$ can be seen formally as an operator from $L_p(\N)$ to $L_p(\N)^{\mathbb{N}}$ by the usual expression : $$Tf(x) = \int_{\mathbb{R}} k(x,y)f(y)dy.$$
We will show in this section that $T$ is bounded from $L_p(\N)$ to $RC_p(\N)$. 

To apply the theorem, we have to include $\M$ and $\ell^2$ in a von Neumann algebra such that their images commute. Let $\Fn_{\infty}$ be the free group with countably many generators $(g_n)_{n\in\mathbb{N}}$, $\lambda$ its left regular representation and $C(\Fn_{\infty})$ its associated von Neumann algebra. Define $\Mt = C(\Fn_\infty) \otimes \M$, $\Nt = L_{\infty}(\mathbb{R}) \otimes \Mt$ and the inclusion maps $i_1 : \M \to \Mt$ and $i_2 : \ell^2 \to \Mt$ by:
\begin{center}
$i_1(x) = 1 \otimes x$ and $i_2(h) = \Sum{k=0}{\infty} h_k \lambda(g_k) \otimes 1$.
\end{center}
From now on $\M$ will also designate $1 \otimes \M$ and can be considered as a von Neumann subalgebra of $\Mt$ since $\tau_{\Mt}$ coincides with $\tau_{\M}$ on $\M$. Note that the image of $i_2$ is in $\M' \cap \Mt$. Let $\kt : (x,y) \mapsto i_2(k(x,y))$ be the kernel of a $\CZ$ operator $\Tt : L_p(\N) \to L_p(\Nt)$. The image of $i_2$ is in $\M'$. The following lemma (Khintchine's inequalities for the free group) is the crucial point of the construction (see \cite{Pisier1},\cite{Pisier2}).

\begin{lemma} \label{ineq:khintchine}
Let $a = (a_n)_{n\in\mathbb{N}}$ be a sequence in $L_p(\M)^{\mathbb{N}}$, then for all $1 \leq p \leq \infty$: $$\Norm{a}_{\RC{p}{\M}} \approx \Norm{\Sum{n = 0}{\infty} \lambda(g_n) \otimes a_n}_{L_p(\Mt)}$$
\end{lemma}

\begin{corollary} \label{cor:TTt}
Let $T$ be a $\CZ$ operator associated with a kernel $k$ taking values in $\ell^2$ an $\Tt$ as defined above, then for all $1\leq p\leq \infty$ and all $f \in L_p(\N)$: $$ \Norm{Tf}_{\RC{p}{\N}} \approx \Norm{\Tt f}_{L_p(\Nt)}.$$
\end{corollary}

\emph{Proof}. Notice that if $Tf = (T_if)_{i\in\mathbb{N}}$ in $\RC{p}{\N}$ then $\Tt f =  \Sum{i = 0}{\infty} \lambda(g_n) \otimes T_if$ in $L_p(\Nt)$ and apply the previous lemma.
\hfill{$\qed$}

\begin{property}
Let $T$ be a $\CZ$ operator associated to a kernel $k$ taking values in $\ell^2$. Suppose furthermore that $T$ is bounded from $L_2(\N)$ to $\RC{2}{\N}$ and that $k$ satisfies the smoothness and size conditions then the operator $\Tt$ defined above is bounded from $L_1(\N)$ to $L_{1,\infty}(\Nt)$.
\end{property}

\emph{Proof}. We only have to check that Theorem $\ref{th:main}$ can be applied to $\Tt$. Khintchine's inequality for $p = \infty$ imply that for all $h \in \ell^2$, $\Norm{h} \approx \Norm{i_2(h)}_{\Mt}$ so $\kt$ verifies the size and smoothness conditions. Furthermore, corollary $\ref{cor:TTt}$ applied to $p = 2$ gives the boundedness condition on $L_2$. So all the hypothesis are verified and $\Tt$ is bounded from $L_1(\N)$ to $L_{1,\infty}(\Nt)$.
\hfill{$\qed$}

\begin{corollary}
Let $T$ be a $\CZ$ operator associated to a kernel $k$ taking values in $\ell^2$. Suppose furthermore that $T$ is bounded from $L_2(\N)$ to $\RC{2}{\N}$ and that $k$ satisfies the smoothness and size condition. Then $T$ is bounded from $L_p(\N)$ to $\RC{p}{\N}$ for all $1<p\leq 2$.
\end{corollary}

\emph{Proof}. From the previous property, it is clear that $\Tt$ is bounded from $L_p(\N)$ to $L_p(\Nt)$ by real interpolation. Which is enough to conclude thanks to corollary $\ref{cor:TTt}$.
\hfill{$\qed$}

\begin{remark}
$T$ is not bounded from $L_p(\N)$ to $C_p(\N)$, for $1 < p < 2$. It is necessary to consider $RC_p(\N)$.
\end{remark}

\emph{Proof}. We will construct a counter-example thanks to the Littlewood-Paley kernel (the same general idea can again be found in \cite{czoperator}). 

Let $\psi$ be a $C^{\infty}$ function supported in $(1,2)$, bounded by $1$ and constant equal to one on $(5/4,7/4)$. For any $i \in \mathbb{N}^*$, let $\phi_i : t \mapsto \psi(it)$. Consider the kernel $k : \mathbb{R}^2 \to \ell^2$ defined by $k_i(x,y) = \phih_i(x - y)$ where $\phih_i$ denotes the Fourrier transform of $\phi_i$. It is standard that $k$ verifies the smoothness estimate for $\gamma = 1$ and the size condition. The $\CZ$ operator $T$ associated to $k$ is bounded on $L_2$ by Plancherel's theorem.

Let $\M = \B(\ell^2)$. As always, let $\N = L_{\infty}(\mathbb{R}) \otimes \M$. For all $p \geq 1$, $T$ can be seen as an operator from $L_p(\N)$ to $C_p(\N)$. For all $k > 0$, let $g_k$ be the Fourier transform of  $\I_{(\frac{5}{4}2^{k-1},\frac{5}{4}2^{k-1} + 1/2)}$. Notice that $\phih_i \ast g_k = \delta_{}$.
Fix $m > 0$ an integer. Let $f_m = \Sum{k=1}{m} g_k \otimes e_{1,k}$ then $Tf_m = (g_k \otimes e_{1,k})_{e_k\in\mathbb{N}^*}$. It results that $\Norm{f_m}_{L_p(\N)} = m^{1/2}\Norm{g_1}_p$ and $\Norm{Tf_m}_{C_p(\N)} = m^{1/p}\Norm{g_1}_p$. So $T$ is not bounded for $p<2$.
\hfill{$\qed$}

\section{$L_p$-pseudo-localisation}

In this section, we will recover in this context the main result from \cite{Hytonen} which is an $L_p$-version of pseudo-localisation. The proof heavily relies on estimates and definitions introduced in section $2.$ Precisely, we will show the following theorem, using notations from theorem \ref{th:pseudoloc}.

\begin{theorem} \label{th:Lp pseudoloc}
Let $T$ be a $\CZ$ operator associated with a kernel $k$ with Lipschitz parameter $\gamma$, verifying the size condition, taking values in $\M' \cap \Mt$ and bounded from $L_2(\N)$ to $L_2(\Nt)$. Then for all $1<p<\infty$, $s \in \mathbb{N}$ and $f \in L_p(\N)$, there exists $\theta_p > 0$ such that : 
\begin{center}
$\Norm{A_{f,s}^{\bot}(Tf)}_p \les 2^{-\gamma s \theta_p }\Norm{f}_p$ and $\Norm{(Tf)B_{f,s}^{\bot}}_p \les 2^{-\gamma s \theta_p }\Norm{f}_p,$
\end{center}
where the implied constant depends on $p$.
\end{theorem}

We will use the decomposition we obtained in the $L_2$-case in (\ref{Decomposition}):
\[
A_{f,s}^{\bot}Tf = A_{f,s}^{\bot}\big(\Sum{i=0}{\infty} \Phi_i f + \Sum{i=1}{\infty}\Psi_i f\big).
\]

The crucial result, proved in the next sections, is the following :

\begin{property}\label{es:crucial}
For all $p\in (1,\infty)$ the sequences of operators $(\Phi_i)_{i\in\Nb}$ and $(\Psi_i)_{i\in\Nb}$ are uniformly bounded on $L_p(\N)$.
\end{property}

\emph{Proof\ of\ the\ theorem.} Fix $p\in (1,\infty)$. By complex interpolation, there exist $q\in (1,\infty)$ and $\theta_p > 0$ such that $\Norm{\Phi_i}_{\B(L_p(\N))} \leq \Norm{\Phi_i}_{\B(L_2(\N))}^{2\theta_p}\Norm{\Phi_i}_{\B(L_q(\N))}^{1-2\theta_p}.$ 

Now, it suffices to plug in the previous proposition and proposition \ref{es:phi_i psi_i} to get the expected estimate. We obtain: $\Norm{\Phi_i}_{\B(L_p(\N))} \les 2^{-\theta_p(i+s)}$. A similar estimate is true for $\Psi_i$. This is enough to conclude thanks to the decomposition (\ref{Decomposition}) mentionned above.
\hfill{$\qed$}

\subsection{$L_p$-boundedness of $\Phi_i$}

\begin{lemma}
Consider a $\CZ$ operator $T$ with kernel $k$, and an integer $m$. Then, for any $p \in (1,\infty)$, the kernel $k' : (x,y) \mapsto \I_{x\in 5Q_{y,m}}k(x,y)$ defines a bounded operator $T'$ on $L_p(\N)$.
\end{lemma}

\emph{Proof}. We introduce a sequence of Rademacher variables indexed by the set of dyadic cubes of size $2^{-m}$ to express the kernel $k'$ in a suitable way. Let $\Omega$ be the probability space $\{-1,1\}^{\Q_m}$ where $\Q_m$ is the set of dyadic cubes of size $2^{-m}$ and $\Nt = \N \overline{\otimes} L_{\infty}(\Omega)$. Consider $\N$ as a subalgebra of $\Nt$. We will take a functional approach and work in $\Nt$ from now on. Denote by $\e_Q$ the $Q$-coordinate in $\Omega$. Define $g_m(x,\omega) = \Sum{Q\in\Q_m}{} \I_Q(x)\e_{Q}(\omega)$, $l_m(x,\omega) = \Sum{Q\in\Q_m}{} \I_{5Q}(x)\e_{Q}(\omega)$. Let $G_m$ (resp. $L_m$) be the multiplication by $g_m$ (resp. $l_m$) and $\E$ be the conditionnal expectation onto $L_p(\N)$. Note that $\Norm{G_m}_{\B(L_p(\N))} = \Norm{g_m}_{\infty} = 1$ and $\Norm{L_m}_{\B(L_p(\N))} = \Norm{l_m}_{\infty} = 5^n$.

It is easily checked that $T' = \E L_mTG_m$ since $$k'(x,y) = \Int{\Omega}{} l_m(x,\omega)k(x,y)g_m(y,\omega) d\omega. $$

So $T'$ can be extended to a bounded operator on $L_p(\N)$.
\hfill{$\qed$}

\begin{property}
Let $p\in (1,\infty)$. The operators $\Phi_i$ are uniformly bounded on $L_p(\N)$.
\end{property}

\emph{Proof}.
We identify $L_p(\N)$ quasi-isometrically with a subspace of $RC_p(\Nt)$, thanks to Burkholder-Gundy inequality (theorem $2.1$ in \cite{PisierXu2}):  $\Norm{f}_p \approx \Norm{df}_{RC_p(\Nt)}$.
Recall that $\Phi_i = \Sum{k\in\Z}{} \Delta_{k-s+i}T_{k-s}\Delta_{k} = \Sum{k\in\Z}{} \Delta_{k-s+i}\big(T - \E L_{k-s}TG_{k-s}\big)\Delta_{k}$, where the $L_k$ and $G_k$ are given by the lemma above. We can write $\Phi_i = \Delta T - \Delta\E LTG$ where $Gx = (G_{k-s}x_k)_{k\in\Z}$, $Lx = (L_{k-s}x_k)_{k\in\Z}$, $\Delta x = (\Delta_{k-s+i}x_k)_{k\in\Z}$ and $Tx = (Tx_k)_{k\in\Z}$. $L$ and $G$ are still bounded multiplications by unconditionnality of $RC_p(\Nt)$. $T$ is completely bounded on $RC_p(\Nt)$ as a diagonal $\CZ$ operator, $\E$ as a conditional expectation from $RC_p(\Nt)$ to $RC_p(\N)$ and $\Delta$ by Stein's inequality ( Theorem $2.3$ in \cite{PisierXu2}):
\[
\Norm{(\E_n(x_n))}_{L_p(\Nt,\ell_2^c)}\leq \Norm{(x_n)}_{L_p(\Nt,\ell_2^c)} \text{ and } \Norm{(\E_n(x_n))}_{L_p(\Nt,\ell_2^r)}\leq \Norm{(x_n)}_{L_p(\Nt,\ell_2^r)}.
\]
Hence, the $\Phi_i$ are uniformly bounded.

\hfill{$\qed$}

\subsection{$L_p$-boundedness of $\Psi_i$}

Recall that,
\[
\Psi_i = \Sum{k\in \Z}{} \Delta_{k-i}S_{k,s}\Delta_{k+s} = \Sum{k\in \Z}{} \Delta_{k-i}(T_{k,s} + R_{k,s})\Delta_{k+s} = A_i + B_i
\]
where $A_i = \Sum{k\in \Z}{} \Delta_{k-i}T_{k,s}\Delta_{k+s} = \Sum{k\in \Z}{} \Delta_{k-i}T_k\Delta_{k+s}$ (by \ref{eq:Tk'}) and $B_i = \Sum{k\in \Z}{} \Delta_{k-i}R_{k,s}\Delta_{k+s}$.
The previous section tells us that $A_i$ is bounded from $L_p(\N)$ to $L_p(\N)$, with a control on the norm independant from $i$. So what is left to prove is the following proposition.

\begin{property}
Let $p \in (1,\infty)$, then $B_i$ is bounded from $L_p(\N)$ to $L_p(\N)$, uniformly in $i$.
\end{property}

\emph{Proof}. Using notations from the previous section, we will again consider $B_i$ as a partial operator from $RC_p(\Nt)$ to $RC_p(\Nt)$. Let $\phi$ be a Schwarz function from $\R^{n}$ to $\R$ such that $\I_{\A_0}(x,y) \leq \phi(x-y) \leq \I_{\md{x-y}>1/2}$ for all $x,y \in \R^n$. Let $U$ be a smooth function from $\Sb^{n-1}$ to the unit circle of $\Cb$ with average $0$. Such $U$ exists, take for example $U_{\lambda}(x) = e^{2i\arctan(\lambda x_1)}$. Since $\arctan$ is an odd function, the average of $U_{\lambda}$ is real. Now, note that $U_{\lambda}$ goes to $-1$ when $\lambda$ goes to $\infty$ and to $1$ when $\lambda = 0$ so there must exist a suitable $U_{\lambda}$. Extend $U$ to $\R^n$ by $U(0) = 0$ and $U(x) = U(x/\md{x})$ otherwise. Write,

\begin{align*}
r_{k,s}(x,y) &= \I_{\A_k}(x,y)\dfrac{K(x)2^{-(k+s)\gamma}}{I_{n,\gamma}\md{x-y}^{n+\gamma}} \\
&= \dfrac{K(x)2^{-s\gamma}}{I_{n,\gamma}}.\I_{\A_k}(x,y).\dfrac{\overline{U}(x-y)\phi(2^k(x-y))}{\md{2^k(x-y)}^{\gamma}}.\dfrac{U(x-y)}{\md{x-y}^n}
\end{align*}

Denote by $M_K$ the multiplication by $\dfrac{K(.)2^{-s\gamma}}{I_{n,\gamma}}$, $M_K$ is bounded (see section $2.3$). Define:
\[
F(x) = \dfrac{\overline{U}(x)\phi(x)}{\md{x}^{\gamma}}.
\]
Since $\phi$ is zero in a neighbourhood of the origin, $F$ has no singularity and is a Schwarz function. So $F$ is the Fourier transform of an $L_1$ function $\widehat{F}$. Note also that $\I_{\A_k} = \I_{x\in 5Q_{y,k}} - \I_{x\in 3Q_{y,k}}$, we will only prove the boundedness replacing $\I_{\A_k}$ by $\I_{x\in 5Q_{y,k}}$ since $\I_{x\in 3Q_{y,k}}$ is similar, denote the associated operator $B_i'$. Define also:
\[
h(x) = \dfrac{U(x)}{\md{x}^n}
\]
and $H$ the convolution operator associated to $h$. Since $U$ has $0$ average, $H$ is bounded from $RC_p(\Nt)$ to $RC_p(\Nt)$ for any $p\in (1,\infty)$. Indeed, $H$ can then be written as an average of directionnal Hilbert operators (see \cite{Duoandikoetxea}). Using notations from the previous proposition and lemma:
\begin{align*}
\I_{x\in 5Q_{y,k}}F(2^{k}(x-y))h(x-y) &= \Int{\Omega}{} l_k(x,\omega)F(2^{k}(x-y))h(x-y) g_k(y,\omega)d\omega \\
& =\Int{\Omega}{}\Int{\R^n}{} \widehat{F}(t)l_k(x,\omega)e^{2^kitx}h(x-y)e^{-2^kity}g_k(y,\omega) dtd\omega.
\end{align*}

We will now translate this equality in terms of operators. Take again $L$ and $G$ from the previous part, the diagonal operators associated to the multiplication by $(l_k)_{k\in\Z}$ and $(g_k)_{k\in\Z}$, $M_t$ defined by $M_t x = (e^{2^ki t.}x_k)_{k\in\Z}$, $\E$ the expectation from $RC_p(\Nt)$ to $RC_p(\N)$, and $\Delta$ such that $\Delta x = (\Delta_{k-s-i}x_k)_{k\in\Z}$. Then,

\[
B'_i = M_K\Int{\R^n}{} \widehat{F}(t)\Delta\E L M_t H M_{-t} G dt.
\]
As we have seen in the previous proof, $M_t$,$M_{-t}$,$L$ and $G$ are bounded by unconditionnality of $RC_p(\Nt)$. More precisely, $\Norm{M_t}_{\B(RC_p(\Nt))} = \Norm{M_{-t}}_{\B(RC_p(\Nt))} = \Norm{G}_{\B(RC_p(\Nt))} = 1$ and $\Norm{L}_{\B(RC_p(\Nt))} = 5^n.$ $\Delta$ is bounded thanks to Stein's inequality and we have have constructed $F$ and $H$ such that $\widehat{F}$ is $L_1$ and $H$ is bounded on $RC_p(\Nt)$. Consequently,
\[
\Norm{B_i}_{\B(L_p(\N))} \lesssim 2\Norm{\Delta}_{\B(RC_p(\N))} 5^n\Norm{\widehat{F}}_1\Norm{H}_{\B(RC_p(\N))}.
\]

\hfill{$\qed$}

\nocite{MeiParcet}

\bibliographystyle{apalike}
\bibliography{bibli}

\emph{Laboratoire de mathématiques Nicolas Oresme, Université de Caen Normandie, 14032 Caen Cedex, France.}

\emph{E-mail address}: leonard.cadilhac@unicaen.fr

\end{document}